\documentclass[10pt,oneside,leqno]{amsart}
\usepackage{amsxtra}
\usepackage{amsopn}
\usepackage{amsmath,amsthm,amssymb}
\usepackage{amscd}
\usepackage{amsfonts}
\usepackage{latexsym}
\usepackage{verbatim}
\usepackage{pb-diagram}
\usepackage{diagrams}

\theoremstyle{plain}
\newtheorem{theorem}{Theorem}[section]
\newtheorem*{theorem*}{Theorem}
\newtheorem{definition}[theorem]{Definition}

\newtheorem{prop}[theorem]{Proposition}
\newtheorem{cor}[theorem]{Corollary}
\newtheorem{rem}[theorem]{Remark}
\newtheorem{ex}[theorem]{Example}
\newtheorem*{mt*}{Main Theorem}
\newcommand\C{{\mathbb C}}

\newcommand\R{{\mathbb R}}
\newcommand\RR{{\mathbb R}}
\newcommand\Z{{\mathbb Z}}

\newcommand\ad{{\rm ad}}

\newcommand{\cg}{{\mathfrak{g}}}

\begin{document}

\title[On the de Rham cohomology of solvmanifolds]{On the de Rham cohomology of solvmanifolds}
\author{Sergio console and Anna Fino}
\date{\today}
\address{Dipartimento di Matematica G. Peano \\ Universit\`a di Torino\\
Via Carlo Alberto 10\\
10123 Torino\\ Italy} \email{sergio.console@unito.it}
 \email{annamaria.fino@unito.it}
\thanks{This work was supported by the Projects MIUR ``Riemannian Metrics and Differentiable Manifolds'',
``Geometric Properties of Real and Complex Manifolds'' and by GNSAGA
of INdAM}
\begin{abstract} By using  results by D. Witte \cite{Witte}  on the superigidity of lattices in solvable Lie groups we get a different proof of a recent  remarkable result obtained by D. Guan \cite{Guan}  on the de Rham cohomology of  a compact solvmanifold, i.e. of a quotient of a connected and simply connected solvable Lie group $G$ by a lattice $\Gamma$. This result can be applied   to compute the Betti numbers of a compact solvmanifold $G/\Gamma$  even in the case that   the solvable Lie group $G$ and the lattice $\Gamma$ do  not satisfy the Mostow condition.\end{abstract}
\maketitle
\section{Introduction}

Let $M$ be a compact solvmanifold, i.e. a quotient of a connected and simply connected solvable Lie group $G$ by a lattice $\Gamma$.

Denote by  $Ad_G (G)$  (respectively,  $Ad_G (\Gamma)$) the subgroup of $GL ({\mathfrak g})$ generated by $e^{{\mbox {ad}}_X}$, for all $X$ in  the Lie algebra   ${\mathfrak g}$ of $G$  (respectively,  in the Lie algebra of $\Gamma$).

 It is well known  that if $G$ is  a simply connected  solvable Lie group,  then  $Ad_G (G)$ is a solvable algebraic group and ${\mbox {Aut}} (G) \cong {\mbox {Aut}} ({\mathfrak g})$. 
We will denote by ${\mathcal A} (Ad_G (G))$  and ${\mathcal A} (Ad_G (\Gamma))$ the real algebraic closures of $Ad_G (G)$ and $Ad_G (\Gamma)$ respectively.

In general, as a consequence of  the Borel density theorem (see \cite[Corollary 4.2]{Witte2} and Theorem~\ref{Borel} here) applied to the adjoint representation, one has that if $\Gamma$ is a lattice of a connected solvable Lie group $G$, then ${\mathcal A} (Ad_G (G)) = T_{cpt} {\mathcal A} (Ad_G (\Gamma))$,  where $T_{cpt}$ is any maximal compact torus  of  ${\mathcal A} (Ad_G (G))$.

If ${\mathcal A} (Ad_G (G)) = {\mathcal A} (Ad_G (\Gamma))$, i.e.,  if  $G$ and $\Gamma$ satisfy the {\em Mostow condition},   then  the de Rham  cohomology $H^*_{dR} (M)$ of  the compact  solvmanifold $M = G / \Gamma$  can be computed  by the  Chevalley-Eilenberg cohomology $H^* ({\mathfrak g})$ of the Lie algebra ${\mathfrak g}$ of $G$ (see \cite{Mostow2} and \cite[Corollary 7.29]{Raghunathan}); indeed, one has the isomorphism
\begin{equation} \label{isocohrag}
H^*_{dR} (M) \cong H^* ({\mathfrak g}).
\end{equation}

In the special case that  $G$ is nilpotent, or equivalently that $M = G / \Gamma$ is a nilmanifold, then the Lie group $Ad_G (G)$ is unipotent and therefore its real algebraic closure  ${\mathcal A} (Ad_G (G))$ coincides with $Ad_G (G)$. Moreover, the real algebraic closure ${\mathcal A} (Ad_G (\Gamma))$  of  the lattice $Ad_G (\Gamma)$ coincides  with  $Ad_G (G)$. For nilmanifolds the isomorphism \eqref{isocohrag} was previously shown by K. Nomizu \cite{Nomizu} by using  the Leray-Serre spectral sequence.

Nilpotent Lie groups can be viewed as a special class of solvable Lie groups, called the {\em completely solvable} ones. We recall that a solvable Lie group $G$ is completely solvable if  all the linear operators ${\mbox {ad}}_X: {\mathfrak g} \to {\mathfrak g}$ have only real eigenvalues (zero eigenvalues in the nilpotent case),  for all $X \in {\mathfrak g}$. 
Like  in the nilpotent case,  also for  a completely solvable Lie group $G$,   any lattice $\Gamma$ of $G$ satisfies  the Mostow condition, so  \eqref{isocohrag}   holds and the isomorphism  was first proved by A. Hattori \cite{Hattori} by using a different  method.

In the case that ${\mathcal A} (Ad_G (G)) \neq {\mathcal A} (Ad_G (\Gamma))$, it is in general  very difficult to compute  the  de Rham cohomology of $G/\Gamma$.

As far as we know,  the only known result has been obtained recently  by  D. Guan  and it   has been applied to  the classification of compact complex homogeneous manifolds with pseudo-k\"ahler structures (see \cite[Main Theorem 1]{Guan} and \cite{Guan2,Guan3}).  He proved that, if $M = G / \Gamma$ is a compact  solvmanifold,
then there exists a  finite covering space  $\tilde M =  G / \tilde \Gamma$, i.e., $\Gamma/ \tilde \Gamma$ is a  finite group, such that  $\tilde M = \tilde G / \tilde \Gamma$, where 
$\tilde G$ is  another simply connected solvable real Lie group $\tilde G$, 
 diffeomorphic to $G$,  and   satisfiying    the Mostow  condition $${\mathcal A} (Ad_{\tilde G} (\tilde G)) = {\mathcal A} (Ad_{\tilde G} (\tilde \Gamma)). $$

Guan proved  the previous theorem by using the results by V. V. Gorbatsevich \cite{Gorbatsevich, Gorbatsevich2}  on Malcev's splitting  and algebraic splitting  of a solvable Lie group.
The construction of the Malcev's splitting (also called semisimple splitting)   has been introduced by A. Malcev \cite{Maltcev} and later improved by L. Auslander \cite{Auslander}. The  general idea  consists in   embedding   solvable and general simply connected Lie groups   into  split groups, i.e.,  in groups  which are representable as semidirect products $T \ltimes N$, where $N$ is the nilradical and $T$ is an abelian subgroup acting on $N$ by semisimple automorphisms.
 
 The nilradical of the semisimple splitting is called in the literature the {\em nilshadow}.  The traditional nilshadow construction (see also \cite{AT})  kills an entire maximal torus.  D. Witte in  \cite{Witte} refined the previous construction and he introduced a method for   killing  a specific compact subtorus of the maximal torus.

 By using this  Witte's result, we  find  a different proof of  Guan's theorem. More precisely, we show the following:

 \smallskip
 
 {\bf Main Theorem}  {\it  Let $M = G / \Gamma$ be a compact solvmanifold,  quotient of  a simply connected solvable Lie group $G$  by a lattice  $\Gamma$ and let $T_{cpt}$ be a   compact torus such that  $$T_{cpt}  {\mathcal A} (Ad_G (\Gamma)) = {\mathcal A} (Ad_G (G)).$$ Then there exists a  subgroup $\tilde \Gamma$ of finite index in $\Gamma$ and a  simply connected   normal subgroup $\tilde G$ of $T_{cpt} \ltimes G$  such that  $${\mathcal A} (Ad_{\tilde G}  (\tilde \Gamma)) = {\mathcal A} (Ad_{\tilde G} ( \tilde G)).$$ Therefore, 
 $ \tilde G /  \tilde \Gamma$ is diffeomorphic to $G /  \tilde \Gamma$ and $H_{dR}^* (G /  \tilde \Gamma) \cong H^*(\tilde {\mathfrak g})$.} 

\smallskip

Observe that $H_{dR}^* (G /  \Gamma) \cong H_{dR}^* (G /  \tilde \Gamma)^{\Gamma/\tilde \Gamma}$ (the invariants by the action of the finite group $\Gamma/\tilde \Gamma$). 

\smallskip

In Section \ref{secmainth} we give a proof of the previous theorem and  we apply it to the special class of {\em almost abelian Lie groups}.
We recall that a  simply connected Lie group $G$  is called almost abelian if its Lie algebra ${\mathfrak g}$ has an abelian ideal of codimension $1$. 
In \cite{Gor} the Mostow condition was studied for simply connected  almost abelian Lie groups and the first Betti number of a compact quotient of a simply connected almost abelian Lie group by a lattice was already determined  in \cite{Bock}.

\smallskip
\noindent {\em{Acknowledgements}} We would like to thank Diego Conti and Daniel Guan  for useful
comments.

\section{General properties of lattices in solvable Lie groups}

Let $G$  be a real Lie group.  We recall that an element $g \in G$ is {\em unipotent} ($Ad$-unipotent)  if the operator $Ad_g$  on the Lie algebra ${\mathfrak g}$ of $G$ is unipotent, i.e. all of its  eigenvalues are equal to $1$ or equivalently there exists a positive integer $k$ such that $(Ad_g - id)^k =0$. An element $g \in G$ is called {\em semisimple} ($Ad$-semisimple)  if the operator $Ad_g$ on the Lie algebra ${\mathfrak g}$ of $G$ is $\C$-diagonalizable.  A subgroup $H \subset G$ is unipotent if all its elements 
are unipotent.

Given two subgroups $H$ and $H'$ of $G$, the commutator $[H, H']$ is the subgroup of $G$ generated by the elements
$$
\{h, k \} = hk h^{-1}k^{-1}, \quad h \in H, k\in H'.
$$
Given a Lie group $G$ one may consider the descending chain of normal subgroups
$$
G_0 = G \supset G_1 = [G, G]  \supset  \dots  \supset G_{i + 1} = [G_i, G] \supset \dots.
$$
The Lie group $G$ is \emph{nilpotent} if the previous descending chain degenerates, i.e. if  $G_i= \{  e \}$, where $e$ is the identity element, for all $i$ greater or equal to some $k$.

One may consider an other  descending chain of normal subgroups
$$
G_{(0)} = G \supset G_{(1)} = [G, G]  \supset  \dots  \supset G_{(i + 1)} = [G_{(i)}, G_{(i)}] \supset \dots.
$$
and $G$ is \emph{solvable} if the previous chain degenerates. In particular a solvable Lie group is {\em completely solvable} if every eigenvalue $\lambda$ of every operator $Ad_g$, $g \in G$, is real.

The unipotent radical $Rad_u G$  of  a Lie group $G$  is its largest  normal subgroup consisting of unipotent elements. In particular  $Rad_u G$ is a nilpotent subgroup of $G$. We recall that  the {\em nilradical} $N$ of a Lie group $G$ is its largest connected nilpotent normal subgroup. In the general case $Rad_u G \subset N$.   The nilradical $N$ is unique and the quotient $G/N$ is   a simply-connected  abelian Lie group. Thus $G$ is an extension of $N$ by an abelian Lie group $\R^k$, i.e. $G$ satisfies the exact sequence
$$1 \to N \to G \to \R^k \to 1,
$$
for some $k \geq 0.$

If $G$ is solvable,  the commutator $[G, G]$ is nilpotent and therefore  the commutator $[G, G]$  is a subgroup of  $Rad_u G$ and thus it is a subgroup of $N$.

Any compact solvable Lie group is commutative and hence it is a torus. Any simply connected solvable Lie group is homeomorphic to a vector space, although the exponential map ${\mbox {exp}}: {\mathfrak g} \to G$ is not necessarily injective or surjective. Moreover, a connected solvable Lie group $G$ is simply connected if and only if $G$ has no non trivial compact subgroups.

A discrete  subgroup $\Gamma$ of  a solvable Lie group $G$ is a lattice if $G/\Gamma$ has a finite invariant measure. 

By \cite[Theorem 3.1]{Raghunathan}   if $G$  is a solvable Lie group with countably many connected components and $H$ be a closed subgroup, then $G/H$ has a $G$-invariant finite measure if and only if $G/H$ is compact. Therefore, a uniform discrete subgroup of a connected solvable Lie group $G$ is always a lattice. The identity component $\Gamma^0$ of a lattice $\Gamma$ of $G$ belongs to the nilradical $N$ of $G$ and the space $N/ N \cap \Gamma$ is compact. 

For any lattice $\Gamma$ of $G$ we have that $\Gamma \cap N$ is a lattice in  the nilradical  $N$ \cite{Mostow, Raghunathan, Auslander}.

In general if  $G$ is a connected Lie group and $G/\Gamma$ is a compact quotient of $G$ by a uniform discrete subgroup, then the normalizer $N_G (\Gamma^0)$ is uniform in $G$ ($\Gamma^0$ denotes the connected component of the identity).

In the case of  solvable Lie groups there is no simple criterion for the existence of a lattice in a connected and simply connected solvable Lie group. The only criterion is that  the Lie group $G$  has to be unimodular  (see \cite[Lemma 6.2]{Milnor}), i.e. the trace ${\mbox {tr}}  (\ad_X)$ vanishes  for any $X$ in the Lie algebra ${\mathfrak g}$ of $G$.  Since the  Lie algebra ${\mathfrak g}$ is unimodular, then  Poincar\'e duality holds for the (co)-homology of ${\mathfrak g}$ (see for instance \cite[Section 5.11]{GHV}).

\section{Linear algebraic groups}

We  recall  that a subgroup $A$ of $GL(n, \R)$ is a {\em real algebraic group}  or Zariski closed if $A$ is the set of  zeros of some collection of regular functions on $GL(n, \R)$, where for regular functions we mean real valued functions on $GL(n, \R)$ for which there is a polynomial  $p$ in $n^2 + 1$ variables, such that $f (g) = p (g_{ij}, det (g^{-1}))$. Indeed, $GL(n, \R)$ can be viewed as a closed subgroup of $SL(n + 1, \R)$ via the embedding $\rho: GL(n, \R) \to SL(n + 1, \R)$ defined by
$$
\rho (A) = \left ( \begin{array} {ccccc} {}&{}&{}&{}&0\\{}& A &{}&&{}.\\{}& {}&{}&&{}.\\{}& {}&{}&{}&0\\ {}&{}&{}&{}&{}\\ 0&\ldots&{}&0&\frac{1}{\det A}
\end{array}
\right).
$$

As a subset of $\R^{n^2}$, a real algebraic group carries both the natural Euclidean topology 
which turns it into a real Lie group and also the Zariski topology induced from $GL(n, \R)$.

A subgroup $A$ of $GL(n, \R)$  is an {\em almost algebraic group}  or  almost  Zariski closed if there is a real algebraic group $B \subset GL(n, \R)$ such that $B^0 \subset A \subset B$, where $B^0$ is the identity component of $B$ in the Euclidean topology. The difference between being real algebraic and almost real algebraic is  only that $B^0$ has always finite index in $B$.

A real algebraic group $U$ is {\em unipotent} if every element $g$  of $U$ is unipotent, i.e.,  there exists a positive integer $k$ such that $(g - I)^k =0$.  A real algebraic torus $T$ is a torus if $T$ is abelian and Zariski connected and every element  $g$ of $T$ is semisimple, i.e. $g$ is $\C$-diagonalizable.

If $G$ is  a Zariski-connected, solvable real algebraic group, then (cf. \cite{Borel}):
\begin{enumerate}

\item the set $U$ of all  unipotent elements of $G$ is a normal subgroup, which is called the {\em unipotent radical} of $G$;

\item for every maximal torus $T$ of $G$, one has $G = T \ltimes U$;

\item $G/U$ is abelian.

\end{enumerate}

In particular,  a Zariski-connected, solvable real algebraic group  $G$ has a decomposition $$G =  (T_{split} \cdot T_{cpt}) \ltimes U,$$ where 
$$
T_{split} = \R^{\#} \times \R^{\#} \times \ldots \times \R^{\#},
$$
with $\R^{\#}$ a $\R$-split torus (i.e. every element of $\R^{\#}$ is diagonalizable over $\R$) via the embedding
$$
\R^{\#} \to GL(2, \R), \quad  x \mapsto  \left ( \begin{array}{cc} x&0\\ 0 &  \frac{1}{x} \end{array} \right),
$$
$T_{cpt}$ the maximal compact subtorus of $T$ and $T_{split}\cap T_{cpt}$ is finite.

In any connected almost algebraic group $G$ all the maximal compact tori are conjugate. 

The {\em (almost) Zariski closure} ${\mathcal A} (H)$ of a subgroup $H$ of an (almost) algebraic group is the unique smallest (almost) algebraic subgroup that contains $H$.

If $K$ is the Zariski closure of $H$, then ${\mathcal A}  (H) = K^0 H$, i.e.,  ${\mathcal A} (H)$ is the union of all connected components of $K$ that intersect $H$.

For instance, if $H$ is connected, then ${\mathcal A}  (H) = K^0$, so the almost Zariski closure of a connected group is connected.

In general,   given a Lie group $G$, we recall that $Ad_G (G)$ is the subgroup of $GL ({\mathfrak g})$ generated by $e^{\ad_X}$, for all $X \in {\mathfrak g}$. Since $Ad \, {\mbox {exp}}  \, X = e ^{\ad_X}$, we have that $Ad_G (G) = Int({\mathfrak g})$, where $Int({\mathfrak g}) = \ad ({\mathfrak g}).$  Moreover $Ad_G G$ is isomorphic to $G/ C$, where $C$ is the center of $G$.

It turns out that if $G$ is  a simply connected  solvable Lie group then  $Ad_G (G)$ is a solvable algebraic group and ${\mbox {Aut}} (G) \cong {\mbox {Aut}} ({\mathfrak g})$.

If $H$ is a subgroup of a connected Lie group $G$, we will denote by ${\mathcal A} (Ad_G (H))$ the (almost)  Zariski closure of $Ad_G (H)$ in the real  algebraic group ${\mbox {Aut}} ({\mathfrak g})$, where $\mathfrak g$ is the Lie algebra of $G$.

 Since $G$ is connected, by \cite[Proposition 3.8]{Witte}  ${\mathcal A} (Ad_G (G))$ acts on $G$ and centralizes  the center of $G$. If $G$ is in addition simply connected,  ${\mbox {Aut}} (G) \cong {\mbox {Aut}} ({\mathfrak g})$ so ${\mathcal A} (Ad_G (G))$ lifts to a group of automorphisms of $G$.

The set of unipotent elements in ${\mathcal A} (Ad_G (G))$  forms a normal subgroup $U$, the unipotent radical, and if $T$ is any maximal reductive subgroup of ${\mathcal A} (Ad_G (G))$ (a maximal algebraic torus), we have a semidirect product decomposition ${\mathcal A} (Ad_G (G)) =T \ltimes U$.

If $\Gamma$ is a lattice of $G$, one has that ${\mathcal A} (Ad_G (\Gamma))$ contains $U$ (cfr.\cite[p. 45]{Raghunathan}) and so  ${\mathcal A} (Ad_G (\Gamma))$ can be written as $S \ltimes U$ where $S  = {\mathcal A} (Ad_G (\Gamma)) \cap T$.

Note that $S$ leaves   the Lie algebra ${\mathfrak n}$ of the nilradical of $G$ invariant  as well as  the Lie algebra of the identity  connected component of $\Gamma$.

\begin{theorem} (Borel density theorem) \cite[Corollary 3.29]{Witte} \label{Borel} Let $\Gamma$ be a lattice of a  simply connected solvable Lie group $G$, then there exists a  maximal  compact torus $T_{cpt} \subset  {\mathcal A} (Ad_G (G))$, such that 
$$
 {\mathcal A} (Ad_G (G)) = T_{cpt}  {\mathcal A} (Ad_G (\Gamma)).
 $$

\end{theorem}

Since for a completely solvable Lie group $G$  the adjoint representation has only real eigenvalues, as a consequence  (see for instance \cite{Witte}),  we have the following:

\begin{prop} \label{compsolv}  Let $\Gamma$ be a lattice of a  simply connected  completely solvable Lie group $G$,
 then ${\mathcal A} (Ad_G (G))$ has no non trivial  connected compact subgroups and thus as a consequence    ${\mathcal A} (Ad_G (G)) = {\mathcal A} (Ad_G (\Gamma))$.
 \end{prop} 
 
 In particular this is the case when $G$ is nilpotent, since in this case $Ad_G (G)$ is unipotent and thus ${\mathcal A} (Ad_G (G))$ coincides with $Ad_G (G)$.

Note that when  $G$ is  a simply connected nilpotent Lie group, by Birkhoff embedding theorem \cite{Auslander} there exists an embedding a monomorphism $B: G \to  U(n, \R)$, where $U(n, \R) \subset SL (n, \R)$ denotes the group of strictly upper triangular matrices of order $G$ (with all diagonal elements equal to $1$),  which allows to treat $G$ as an algebraic group all of whose connected subgroups are algebraic.

In the special case of a solvable Lie group  $G$  for which  all the eigenvalues of $Ad(g)$, $g \in G$,   belong to the  the unit circle, then  the algebraic  torus  $T$ in ${\mathcal A} (Ad_G (G)) = T \ltimes U$ is compact.

\section{de Rham Cohomology of solvmanifolds}

An important result related to lattices  for  compact solvmanifolds is the following:

\begin{theorem} \cite[Theorem 3.6]{Raghunathan}  \label{diffeo} If $\Gamma_1$ and $\Gamma_2$ are lattices in simply connected solvable Lie groups $G_1$ and $G_2$,  respectively, and $\Gamma_1 \cong \Gamma_2$, then $G_1 / \Gamma_1$ is diffeomorphic to $G_2 / \Gamma_2$.
\end{theorem}

As a consequence  two compact solvmanifolds with isomorphic fundamental groups are diffeomorphic.

\begin{ex} \label{ablattice}  {\rm    Let $G$ be the semidirect product $\R \ltimes_{\varphi} \R^2, $ where
$$
\varphi (t) = \left(  \begin{array}{cc} \cos(2 \pi t)&\sin (2 \pi t)\\ -\sin (2 \pi t)&\cos (2 \pi t) \end{array}
 \right ), \quad t \in \R.
 $$
The element $1 \in \R$ and the standard  integer lattice $\Z^2$ generate a lattice $\Gamma$ isomorphic to $\Z^3$.  By applying previous theorem one has that
$ G / \Gamma$ is diffeomorphic to the torus $T^3$.}
\end{ex}

In general $Ad_G (\Gamma)$ and $Ad_G (G)$ do not have the same Zariski closures. A simple example for which ${\mathcal A} (Ad_G (G)) \neq {\mathcal A} (Ad_G (\Gamma))$ is the following.

\begin{ex}  {\rm 
For the previous Example \ref{ablattice}
we have by \cite{Mostow} that $Ad_G (\Gamma)$ is abelian and then ${\mathcal A} (Ad_G (\Gamma))$ is also abelian, but $Ad_G (G)$ is not abelian. 
Indeed,  ${\mathcal A} (Ad_G (G))$ is locally isomorphic to $G$. Then in this case ${\mathcal A} (Ad_G(\Gamma)) \neq {\mathcal A} (Ad_G  (G))$.

}
\end{ex}

 Let $H^* ({\mathfrak g})$ be the Chevalley-Eilenberg cohomology of  the Lie algebra ${\mathfrak g}$ of $G$. In general, there is a natural injective map  $H^* ({\mathfrak g}) \to H^*_{dR} (M)$  which is an isomorphism in the following cases:

\begin{theorem} \label{Mostcoh}Let $G / \Gamma$ be a  compact solvmanifold. 

\begin{enumerate} \item  \cite{Hattori} If $G$ is completely solvable, then $H_{dR}^* (G / \Gamma) \cong H^* ({\mathfrak g})$.

\item \cite[Corollary 7.29]{Raghunathan}  If $Ad_G (\Gamma)$ and $Ad_G (G)$ have the same Zariski (i.e. real algebraic) closure in $GL ({\mathfrak g})$, then $H^*_{dR}  (G / \Gamma) \cong H^* ({\mathfrak g})$.

\end{enumerate}

\end{theorem} 

For the sake of completness we will  include the proof of (2), generalizing  the one given by Nomizu \cite{Nomizu} in the  case of a nilmanifold.

\smallskip

\noindent {\it  Proof of {\rm (2)}.} 
Let $G_{(k)} = [G_{(k-1)}, G_{(k-1)}]$ and $\Gamma_{(k)} = [\Gamma_{(k-1)}, \Gamma_{(k-1)}]$ be the derived series for $G$ and $\Gamma$. Remark that $G_{(k)}$ is nilpotent for any $k\geq 1$ and that since ${\mathcal A} (Ad_G (G))= {\mathcal A} (Ad_G (\Gamma))$ a general result on lattices in nilpotent Lie groups (\cite[Theorem 2.1]{Raghunathan}) implies that $G_{(k)} /\Gamma_{(k)}$ is compact for any $k$. Hence $G_{(k)} /\Gamma\cap G_{(k)}$ is compact for any $k$.

Let $r$ be the last non-zero term in the derived series of $G$. Namely $G_{(r+1)} =(e)$ and $G_{(r)} =:A \neq (e)$. Observe that $A$ is abelian. So $A/A \cap \Gamma:=T^m$ is a  compact  torus.

 Thus, $\overline M:=G/A\Gamma$ is a compact solvmanifold with dimension smaller than $M:=G/\Gamma$ and   $T^m\hookrightarrow M \stackrel{\pi}\to \overline M$ is a fibration. 

Next, consider the Leray-Serre spectral sequence $E_{*}^{p,q}$ associated with the above fibration. One has 
\[\begin{array}{l} 
E_2^{p,q}=H_{dR}^p(\overline M, H_{dR}^q(T^m))\cong H_{dR}^p(\overline M) \otimes \bigwedge^q \RR^m\, , \\
E_{\infty}^{p, q}\Rightarrow
H^{p+q}_{dR}(M)\, .
\end{array}\]

The idea is to construct a second spectral sequence $\tilde E_{*}^{p,q}$ which is the Leray-Serre spectral sequence relative to the complex of $G$-invariant forms (i.e., the Chevalley-Eilenberg complex). Since $\bigwedge^* \cg^*$ is subcomplex of   $\bigwedge^* M$, $\tilde E_{*}^{p,q} \subseteq \ E_{*}^{p,q}$ and 
\[\begin{array}{l} 
\tilde E_2^{p,q}=H^p(\cg/{\mathfrak a}) \otimes {\bigwedge}^q \RR^m\, , \\
\tilde E_{\infty}^{p, q}\Rightarrow
H^{p+q} (\cg)\, ,
\end{array}\]
where $\mathfrak a$ denotes the Lie algebra of $A$.

Now, since $\overline M$ is a  compact solvmanifold of lower dimension than $M$, an inductive argument on the dimension shows that $H_{dR}^p(\overline M)\cong H^p(\cg/{\mathfrak a})$ for any $p$. Thus $E_2=\tilde E_2$ and one has the equality $E_\infty = \tilde E_\infty$. Therefore, $H_{dR}^k
(M)=H_{dR}^k (G/\Gamma)
\cong H^k ({\mathfrak g})$ for any $k$. 
\hfill{$\Box$}
\smallskip

\noindent \begin{rem}  {\rm  By Proposition \ref{compsolv} it turns out that (1) is then a particular case of  (2).

Moreover, one can have the isomorphism $H^*({\mathfrak g}) \cong H^*_{dR} (G /\Gamma)$ even if $Ad_G (\Gamma)$ and $Ad_G (G)$ do not have the same algebraic closure  (see Example \ref{hyperell}).} \end{rem}

For a solvable Lie algebra $\mathfrak g$ one has always  that $[{\mathfrak g}, {\mathfrak g}] \neq \mathfrak g$. Thus,  as a consequence,  the  first Betti number $b_1(G/\Gamma)$ of a compact solvmanifold $G/\Gamma$ satisfies the inequality:
$$
b_1(G/\Gamma) \geq 1.
$$
(see also \cite[Corollary 3.12]{Bock}). In the particular case of a nilmanifold $G/\Gamma$, one has $b_1(G/\Gamma) \geq 2.$

\section{Semisimple   splitting and nilshadow}

Let $G$ be a  connected and simply-connected solvable Lie group.

 We recall that $G$ can be viewed as  the extension of the nilradical  $N$ by an abelian Lie group $\R^k$, but in general this extension  does not split.

If $G$ is a solvable algebraic  real group, then  $G = A  \ltimes U$, where $U$ is a unipotent radical and $A$ is an abelian subgroup, consisting of semisimple elements  \cite{Merzilyakov}. Thus in this case $G$ splits.

\begin{definition} A solvable simply connected Lie group $G$ is {\em splittable} if  $G$ can be represented in  the form of a semidirect product $G  = A \ltimes_{\varphi} N$, where $N$ is the nilradical of $G$, $A$ is an abelian subgroup and the subgroup $\varphi (A) \subset  {\mbox {Aut}} (N)$ consists of semisimple automorphisms.
\end{definition}

By definition $\ker \varphi$ is a discrete subgroup of $T$.

Malcev introduced for an arbitrary simply connected Lie group the general  notion of {\em splitting}.

\begin{definition} A  {\em semisimple splitting}   or a {\em Malcev splitting} of a simply connected solvable Lie group $G$ is an imbedding $\iota: G \to  M (G)$ of the group $G$ in a splittable simply-connected   solvable Lie group $M (G)$, such that if $M (G) = T_G \ltimes U_G$, where $U_G$ is the nilradical of $M (G)$ and $T_G$ is an abelian subgroup acting on $U_G$ by semisimple automorphisms, then
\begin{enumerate}

\item $M (G) =\iota (G) U_G$ (product of subgroups);

\item $M (G) = T_G \ltimes \iota (G)$ (semidirect product of $T_G$ with the normal subgroup $\iota (G)$).

\end{enumerate}

\end{definition}

\medskip

Gorbatsevich proved in \cite{Gorbatsevich} that, for any simply connected  solvable Lie group, a  Malcev splitting  exists and it is unique. 

Then the  general idea for  semisimple splitting  is embedding solvable and general simply connected Lie groups in splits groups.

We recall  briefly the construction. 

Let $G$ be a simply connected solvable Lie group and $Ad: G \to {\mbox {Aut}} ({\mathfrak g})$ the adjoint representation of $G$ in the automorphism group ${\mbox {Aut}}  ({\mathfrak g})$ of the Lie algebra  $\mathfrak g$ of $G$.  The real algebraic group ${\mbox {Aut}} ({\mathfrak g})$ splits as follows:
$$
{\mbox {Aut}}  ({\mathfrak g}) =(P  T) \ltimes U,
$$
where $P$ is a maximal semisimple subgroup, $T$ is an abelian subgroup of semisimple automorphisms of ${\mathfrak g}$  that commute with $P$  (and the intersection $P \cap T $ is discrete) and $U$ is the unipotent radical of ${\mbox {Aut}}  ({\mathfrak g})$. 

The group $Ad_G (G)$ is a connected solvable normal subgroup of ${\mbox {Aut}} ({\mathfrak g})$ and thus it is contained in the radical $T \ltimes U$ of the algebraic group ${\mbox {Aut}} ({\mathfrak g})$.

Let $\theta: T \ltimes U  \to T$ be the epimorphism induced by projection and $T_G = \theta (Ad_G (G))$. Thus $T_G$ is a connected abelian group of semisimple automorphisms of ${\mathfrak g}$ and therefore of $G$ itself.

The connected solvable Lie group $M (G) = T_G \ltimes G$ is precisely  the semisimple splitting  of $G$.

It is possible to show that (see \cite{AB})  if $U_G$ is the nilradical of $M (G)$, then $M (G)$ is the semidirect product $M (G) = T_G \ltimes_{\nu_s} U_G$. Moreover, $M (G) =  G \, U_G$, i.e. $G$ and $U_G$ generate $M(G)$. The new Lie group   $M (G)$ is splittable, since the exact sequence
$$
1 \to U_G \to M(G) \to \R^k \to 1
$$
splits. The nilradical $U_G$ is also called the {\em nilshadow} of $G$. If we denote by
$$
p: M (G) = T_G \ltimes U_G \to T_G, \quad \pi: M (G) = T_G \ltimes U_G \to U_G,
$$
the projections over $T_G$ and $U_G$, then the restriction $p: G \to U_G$ is a group epimorphism and the restriction $\pi: G \to U_G$ is a diffeomorphism.

Moreover, 
$$
[U_G, U_G]  \subset N, \quad [T_G, U_G] \subset N
$$
(see for instance \cite{Starkov}).

Following \cite{Bock} we give an explicit expression of $U_G$ and  of the action $\nu_s$ of $T_G$ on $U_G$. More precisely, if $N$ is the  nilradical of $G$, we have that there exists a maximal connected and simply connected nilpotent Lie subgroup $K$ of $G$ such that $G$ is the  product $K  \cdot   N$ \cite[Proposition 3.3]{Dekimpe}. The nilpotent Lie group $K$ is also called a {\em supplement} for $N$ in $G$ and it is unique up to conjugation. 

There is a natural action of $K$ on $N$ induced  by conjugation in $G$. 
Consider the action $$\tilde \phi: K \longrightarrow {\mbox {Aut }}(G),$$ given by
$$
\tilde \phi (a) (k \cdot n) = k \cdot (I_a \vert_N)_s (n),
$$
where $(I_a \vert_N)_s (n)$ is the semisimple part of the inner automorphism of $N$, obtained by conjugating the elements of $N$  by $a$.

Then, there is an induced action $\phi$ of 
$$
T_G = K / (N \cap K) \cong (KN)/N \cong G/N \cong \R^k.
$$
on $G$ which makes the following diagram commutative: 

\[
\begin{diagram}
\node{K}
\arrow{se,l}{\pi}  \arrow[2]{es,t,2}{\tilde \phi}\node[2]{{\mbox {Aut}}(G)}
\\
\node[2]{T_G= K /(K \cap N)}\arrow{ne,b}{\phi}
\end{diagram}
\]

The semisimple splitting $M(G)$ coincides with the semidirect product
$$
T_G  \ltimes_{\phi} G.
$$
 Indeed,  
one has  an imbedding $\iota: G \longrightarrow T_G  \ltimes_{\phi} G$ of $G$ into  $T_G  \ltimes_{\phi} G$ such that $\iota (G)$ is a closed normal subgroup of $T_G  \ltimes_{\phi} G$.

Moreover, by \cite{Bock} the nilradical of $T_G  \ltimes_{\phi} G$ is:
$$
\begin{array}{lcl}
U_G &=&  \{ \pi (k^{- 1}) \cdot k \, \mid \, k \in K \} \cdot N\\[4pt]
&=&  \{ (\pi (k^{-1}), k \cdot n)  \mid  k \in K, n \in N \} \subset T_G  \ltimes_{\phi} G
\end{array}
$$
 and $T_G  \cap U_G = \{ e \}$. Then  
$$ T_G  \ltimes_{\phi} G = T_G  \ltimes_{\nu_s}  U_G,
$$
coincides with  the  semisimple splitting $M (G)$  of $G$,
where $U_G$ is the nilradical of $M(G)$ and $\nu_s$ is the semisimple automorphism:
$$\nu_s (t) (\pi(k^{- 1}) \cdot (k \cdot n)) : = t \cdot \pi (k^{-1})  \cdot (k \cdot n) \cdot t^{-1} = \pi (k^{- 1}) \cdot (k \cdot \tilde \phi (k_t) (n)),
$$
for any $t \in T_G, k \in K, n \in N$ and any $k_t \in \pi^{-1} (\{ t \} )$. 

In this way, one has the diagonal embedding:
$$
\Delta : G = K \cdot N  \to M(G) = \R^k \ltimes G, k \cdot n  \mapsto  (\pi (k^{-1}),  k \cdot n)
$$
and the nilradical $U_G$ coincides with the  image  $\Delta (G)$. By construction, $U_G$ is diffeomorphic to $G$ but it has a different product that makes it nilpotent.

More explicitely, at the level of Lie algebras we have that there exists a vector space $V \cong \R^k$ such that ${\mathfrak g} = V \oplus \mathfrak n$ and ${\mbox {ad}}(A)_s (B) =0$,  for any $A, B \in V$, where ${\mbox {ad}} (A)_s$ is the semisimple part of ${\mbox {ad}}(A)$.

By considering $V \cong \R^k$ as abelian Lie algebra, then the Lie algebra of the semisimple splitting $M(G)$  of $G$ is 
$ V \ltimes_{\ad(...)_s} {\mathfrak g}$, i.e.
$$
[(A, X), (B, Y) ] = (0, [X, Y] + {\mbox {ad}} (A)_s (Y) - {\mbox {ad}}(B)_s (X)),
$$ 
with nilradical (the Lie algebra of $U_G$)
$$
{\mathfrak u}_G = \{ (- X_V, X), \mid  X \in {\mathfrak g} \},
$$
where $X_V$ denotes the component of $X$ in $V$. Note that ${\mathfrak u}_G$ is isomorphic to ${\mathfrak g}$ only as a vector space.

\begin{ex} {\rm If $G$ is nilpotent, then by definition $M(G) =G$ and if $G$ is splittable, i.e. $G = A \ltimes N$, where $N$ is the nilradical of $G$ and the abelian group $A$ acts on $N$ by semisimple automorphisms. Let $*: A \to {\mbox {Aut}}(N)$ be the natural homomorphism, then the image $* (A) = A^*$ is an abelian subgroup of ${\mbox {Aut}}(N)$ consisting of semisimple elements and $A^*$ can be considered as a subgroup of ${\mbox {Aut}}$.

Therefore  $T_G = A^*$, $M (G) = A^*  \ltimes G$ and the nilshadow $U_G$ is $$
\{ ((a^*)^{-1}, a) \, \mid a \in  A \} \cdot N.
$$
In this case one says that the nilshadow is the image of the \lq \lq diagonal" map
$$
A \ltimes N \to A^* \ltimes  (A \ltimes N), (a, n) \mapsto ((a^*)^{-1}, (a,n)).
$$
}
\end{ex}

By construction there is a strict relation between the semisimple  splitting of $G$ and the algebraic closure 
$ {\mathcal A} (Ad_G (G))$ of $Ad_G (G)$.
Since $G$ is  solvable and simply connected, then ${\mbox {Aut}} (G) \cong {\mbox{Aut}} ({\mathfrak g})$ and  $Ad_G (G)$ is a solvable algebraic group and it is connected as algebraic group. Therefore its closure ${\mathcal A} (Ad_G(G))$, in  $GL ({\mathfrak g})$ has the Chevalley decomposition
\begin{equation} \label{chevalleydec}
{\mathcal A} (Ad_G(G)) = T \ltimes U
\end{equation}
where $U$ is the  unipotent radical (i.e. the largest normal subgroup of ${\mathcal A} (Ad_G(G))$ consisting of unipotent elements) and $T$ is an abelian subgroup  (a maximal algebraic torus) acting by semisimple (i.e. completely reducible)  linear operators. Although $U$ is nilpotent, it is important to note that $U$ need not to be the nilradical of ${\mathcal A} (Ad_G(G))$. Note that $T$ is any maximal reductive subgroup of  ${\mathcal A} (Ad_G (G))$ and it is a standard fact that all such groups $T$ are conjugate by elements of $U$.

The group ${\mbox {Aut}} ({\mathfrak g})$ is algebraic and $Ad_G (G) \subset {\mbox {Aut}} ({\mathfrak g})$. Therefore, the group  $T$ may be considered as a subgroup of ${\mbox {Aut}} ({\mathfrak g}) \cong {\mbox {Aut}} (G)$. There is a natural epimorphism with kernel $U$:
$$
\theta:  T \ltimes U \longrightarrow T.
$$
Let $T^*_G= \theta(Ad_G (G))$ be the image of $Ad_G (G)$ by the previous epimorphism, i.e.  $T^*_G$ can be viewed as the  image of $Ad_G (G)$ in $T$.   

We can  view either $T$ or $T^*_G$  as a group of automorphisms of $G$. By using   the imbeddings
$$
T^*_G \longrightarrow {\mbox {Aut}} (G), \quad T \longrightarrow {\mbox {Aut}} (G)
$$
we may construct  the two semi-direct products
 $$
 T^*_G  \ltimes G \quad {\mbox {and}} \quad   T  \ltimes G.
$$
The semisimple splitting is the universal covering $M (G)= T_G \ltimes G = \R^k \ltimes G$  of the Lie group $T^*_G \ltimes G$.

In general $T_G$ is not an algebraic subgroup of ${\mbox {Aut}} ({\mathfrak g})$ and it could possibly be not closed in the Euclidean topology, but if  $G$ admits lattices, then the subgroup $T_G$ is closed in the Euclidean topology.

\begin{rem} {\rm The semidirect product $T \ltimes G$ splits as a semi-direct product $T \ltimes U_G$, since the nilradical $U_G$  of $M(G)$ coincides with the nilradical of $T \ltimes G$ and  $T \ltimes G$ admits the structure of algebraic group (defined over the real numbers). Moreover, $G$, as a subgroup of $T \ltimes G$, is Zariski dense and hence one may think of $T \ltimes G$  as a kind of  \lq \lq algebraic hull" of $G$.}
\end{rem}

\medskip

\section{The nilshadow map and its applications to the cohomology} \label{secmainth}

Let $G$ be a connected, solvable real algebraic  group, then $G$ is a semidirect product of the form $G = T \ltimes U$, where $T = T_{split} \times T_{cpt}$ is a torus and $U$ is nilpotent.

One can see, following the paper \cite{Witte}, that   the subgroup $T$ acts on $G$ (by conjugation), thus one can construct the semidirect product:
$$
T \ltimes G = T \ltimes (T \ltimes U),
$$
and the nilshadow  ${\Delta} (G)$ can be embedded in $T \ltimes G$ by the map:
$$
(t, u) \to (t^{-1}, t, u) \in T \ltimes (T \ltimes U),
$$
where $t^{-1}$ acts  on $G =T \ltimes U$ by conjugation.

The anti-diagonal imbedding of $T$ sends it into the center of $T \ltimes G$. In this way the nilshadow $U_G$ can be viewed as the image of the map:
$$
\Delta: G \to T \ltimes G,  g \mapsto \Delta(g) = (\pi(g)^{-1}, g),
$$
where $\pi: G \to T$ is the projection, i.e. $U_G = \Delta (G)$.  This antidiagonal embedding of $T$  is the main idea of the nilshadow construction.

For the general case of a solvable (non nilpotent)  Lie group $G$, which may be not real algebraic, one has that 
${\mathcal A} (Ad_G (G)) = T \ltimes U$ is not unipotent, i.e. a maximal torus $T = T_{split} \times T_{cpt}$ of  ${\mathcal A} (Ad_G (G))$ is non trivial. The basic idea  for the construction of  the nilshadow is to kill $T$ in order to obtain  a nilpotent group. In order to do this one can define, following \cite{Witte},  a natural homomorphism $\pi: G \to T$, which is the composition of the homomorphisms:
$$
\pi: G \stackrel{Ad} \to {\mathcal A} (Ad_G (G)) \stackrel{{\mbox{\tiny projection}}}  \longrightarrow T,
$$
 and the map
$$
 \Delta: G \to T \ltimes G,  g \mapsto \Delta(g) = (\pi(g)^{-1}, g).
$$
The traditional nilshadow construction killes the entire maximal torus $T$ in order to get the nilpotent Lie group ${\Delta} (G)$. Witte  introduced in \cite{Witte} a variation, killing only a subtorus $S$ of $T$.
It is well known that, for every subtorus $S$ of a compact torus $T$, there is a torus $S^{\perp}$ complementary to $S$ in $T$, i.e. such that $T = S \times S^{\perp}$.

As a consequence  of \cite[Proposition 8.2]{Witte} now we are able to  prove the Main Theorem.

\smallskip

\noindent {\it Proof  of the Main Theorem.} By \cite[Theorem 6.11, p. 93]{Raghunathan} it is not restrictive to suppose that 
 ${\mathcal A} ( Ad_G (\Gamma))$ is connected. Otherwise we  pass from $\Gamma$ to a finite index subgroup  $\tilde \Gamma$. This  is equivalent to pass from $M = G/ \Gamma$  to the space $ G/ \tilde \Gamma$ which is a finite-sheeted covering of $M$. Let $T_{cpt}$ be a maximal  compact torus   of  ${\mathcal A} (Ad_G G)$  which contains a maximal compact torus $S_{cpt}^{\perp}$ of ${\mathcal A} (Ad_G (\tilde \Gamma))$. There is  a natural
projection from ${\mathcal A} (Ad_G (G))$ to $T_{cpt}$, given by the splitting ${\mathcal A} (Ad_G (G)) = (A \times  T) \ltimes U,$ where $A$
is a maximal $\R$-split torus and  $U$ is the unipotent radical.

 Let $S_{cpt}$ be a subtorus of $T_{cpt}$ complementary to $S_{cpt}^{\perp}$ so that $T_{cpt} = S_{cpt} \times S_{cpt}^{\perp}$.
Let $\sigma$  be the composition of the homomorphisms:
$$
\sigma: G  \stackrel{Ad}  \longrightarrow {\mathcal A} (Ad_G (G))  \stackrel{{\mbox {\tiny projection}}}  \longrightarrow T_{cpt}  \stackrel{{\mbox {\tiny projection}}} \longrightarrow S_{cpt} \stackrel{ x \to x^{-1}}  \longrightarrow S_{cpt}.
$$
One may define the nilshadow map:
$$
\Delta: G \to S_{cpt} \ltimes G, g \mapsto (\sigma (g), g),
$$
which is not a homomorphism (unless $S_{cpt}$ is trivial and then $\sigma$ is trivial), but  one has
$$
\Delta (ab) = \Delta ( \sigma(b^{-1})  a \sigma( b))  \,  \Delta (b),  \quad \forall a,b \in G
$$
and $\Delta (\gamma g) = \gamma \Delta (g)$, for every $\gamma \in \tilde \Gamma, g \in G$.
The nilshadow map $\Delta$ is a diffeomorphism onto it image and then $\Delta (G)$ is simply connected. The product in $\Delta (G)$ is given by:
$$
\Delta (a) \Delta (b) = (\sigma (a), a)  \, (\sigma (b), b) = (\sigma (a) \sigma (b), \sigma (b^{-1})  a \sigma (b) \, b),
$$
for any $a, b \in G$.

By construction  ${\mathcal A} (Ad_G (G))$ projects trivially on $S_{cpt}$ and $\sigma ( \tilde \Gamma) = \{ e \}$. Therefore
$$
\tilde \Gamma = \Delta ( \tilde \Gamma) \subset \Delta (G).
$$
Let $\tilde G = \Delta (G)$.
 By \cite[Proposition  4.10]{Witte}  $S_{cpt}^{\perp}$ is a maximal compact subgroup of ${\mathcal A} (Ad_{\tilde G} (\tilde G))$ and $S^{\perp} \subset {\mathcal A} (Ad_{\tilde G} (\tilde \Gamma))$, therefore  one has ${\mathcal A} (Ad_{\tilde G} (\tilde G)) ={\mathcal A} (Ad_{\tilde G} (\tilde \Gamma))$ as shown in \cite{Witte}.

By using Theorem \ref{diffeo} we have that $G/\tilde \Gamma$ is diffeomorphic to $\tilde G / \tilde \Gamma$. Finally, by applying Theorem \ref{Mostcoh} we have that
$H^*(G / \tilde \Gamma) \cong H^*(\tilde {\mathfrak g})$. By considering the diffeomorphism 
$$
\Delta: G \to \tilde G,
$$
we have that  $\Delta^{-1}$  induces a finite-to-one covering
map  $$\Delta^*: \tilde G/ \tilde \Gamma \to G/ \Gamma.$$
\hfill{$\Box$}

\medskip

\begin{cor} At the level of Lie algebra, if we denote by $X_{\mathfrak s}$ the image  $\sigma_* (X)$, for $X \in {\mathfrak g}$, then the Lie algebra  $\tilde {\mathfrak g}$ of $ \tilde G$ can be identified by
$$
\tilde {\mathfrak g} = \{ (X_{\mathfrak s}, X) \mid X \in {\mathfrak  g} \}
$$
with  Lie bracket:
$$
[( X_{\mathfrak s}, X), (Y_{\mathfrak s}, Y)] = (0, [X, Y] - \ad(X_{\mathfrak s})  (Y) +\ad( Y_{\mathfrak s}) (X)).
$$
\end{cor}

Now we obtain some applications of the Main Theorem, by computing explicitly the Lie group $\tilde G.$

\begin{ex} {\rm {(\bf Nakamura manifold)}  Consider the  simply connected complex solvable  Lie group $G$
 defined by 
$$
G = \left \{ \left(  \begin{array}{cccc} e^z&0&0&w_1\\ 0&e^{-z}&0&w_2\\ 0&0&1&z\\ 0&0&0&1 \end{array} \right) , \,  w_1, w_2, z \in \C \right \}.
$$
The Lie group $G$ is  the semi-direct product $\C  \ltimes_{\varphi} \C^2$, where 
$$
\varphi (z) = \left(  \begin{array}{cc} e^z&0\\ 0&e^{-z} \end{array} \right).
$$
A basis of  complex left-invariant $1$-forms is given by
$$
\phi_1 = dz, \, \phi_2 = e^{- z} d w_1, \, \phi_3 = e^z d w_2
$$
and in terms of the  real basis of left-invariant $1$-forms $(e^1, \ldots, e^6)$ defined by
$$
\phi_1 = e^1 + i  e^2, \, \phi_2 = e^3 + i  e^4, \, \phi_3 = e^5 + i  e^6,
$$
 we obtain the  structure equations:
$$
\left \{ \begin{array} {l}
d e^j = 0, \quad  j = 1,2,\\[4 pt]
d e^3 = - e^{13} + e^{24},\\[4 pt]
d e^4 = - e^{14} - e^{23},\\[4 pt]
d e^5 = e ^{15} - e^{26},\\[4 pt]
d e^6 = e^{16} + e^{25},
\end{array} \right .
$$
where  we denote by $e^{ij}$ the wedge product  $e^i \wedge e^j$.

Let $B \in SL(2, \Z)$ be a unimodular matrix with distinct real eigenvalues: $\lambda, \frac {1} {\lambda}$.
Consider $t_0 = \log \lambda$, i.e. $e^{t_0} = \lambda$.
Then there exists a matrix $P \in GL(2, \R)$ such 
that $$
P B P^{-1} =  \left(  \begin{array}{cc} \lambda &0\\ 0&\lambda^{-1} \end{array} \right).
$$
Let 
$$
\begin{array}{l}
L_{1, 2 \pi} = \Z [t_0, 2 \pi i] = \{ t_0 k + 2 \pi h i ,  h, k \in \Z \},\\ [4pt]
L_2 = \left\{  P \left(  \begin{array}{c}  \mu \\ \alpha \end{array} 
 \right),  \mu,  \alpha \in \Z[i] 
\right \}.
\end{array}
$$
 Then, by \cite{Yamada}  $\Gamma = L_{1, 2 \pi} \ltimes_{\varphi} L_2$  is a lattice of $G$.
 
 Since $G$ has trivial center, we have that $Ad_G (G) \cong G$ and thus it  is a semidirect product $\R^2 \ltimes \R^4$. Moreover,
for  the algebraic closures of $Ad_G (G)$ and $Ad_G (\Gamma)$ we obtain $$
 \begin{array}{l}
 {\mathcal A} (Ad_G G) = (\R^{\#}  \times S^1) \ltimes  \R^4,\\
 {\mathcal A} (Ad_G \Gamma) = \R^{\#}  \ltimes  \R^4,\\
\end{array} 
$$
where  the split torus $\R^{\#}$ corresponds to the action of $e^{\frac{1}{2} (z + \overline z)}$ and the compact torus $S^1$ to the one of $e^{\frac{1}{2} (z - \overline z)}$.

Therefore in this case ${\mathcal A} (Ad_G (G))  = S^1  {\mathcal A} (Ad_G (\Gamma))$ and ${\mathcal A} (Ad_G (\Gamma))$ is connected.
By applying the Main Theorem   there exists a simply connected normal subgroup $\tilde G = \Delta (G)$ of  $S^1 \ltimes G$. The new Lie group $\tilde G$ is obtained  by killing the action of $e^{\frac{1}{2} (z - \overline z)}$. Indeed, 
we get  that 
$$
\tilde G  \cong  \left \{ \left(  \begin{array}{cccc} e^{\frac{1}{2} (z + \overline z)}&0&0&w_1\\ 0&e^{- \frac{1}{2} (z + \overline z)}&0&w_2\\ 0&0&1&z\\ 0&0&0&1 \end{array} \right), w_1, w_2, z \in \C \right \}.
$$

The diffeomorphism between $G/\Gamma$ and $\tilde G/\Gamma$ was already shown in \cite{Yamada}.  Then in this case  one has the isomorphism $H^*_{dR} (G/ \Gamma)  \cong H^* (\tilde {\mathfrak g})$, where $\tilde {\mathfrak g}$ denotes the Lie algebra of $\tilde G$ and the de Rham cohomology of the Nakamura manifold  $G/\Gamma$ is not isomorphic to $H^* ({\mathfrak g})$ (see also \cite{debaT}).}

\end{ex}

\medskip

\begin{ex}  {\rm Let consider the $3$-dimensional solvable Lie group $\R \ltimes \R^2$ with structure equations
 $$
 \left \{  \begin{array} {l}
d e^1 =0,\\
d e^2 = 2 \pi e^{13},\\
d e^3 = - 2 \pi e^{1 2}.
\end{array} \right .
 $$
 is a non-completely solvable Lie group which admits a compact quotient and the uniform discrete subgroup is of the form
  $\Gamma =  \Z \ltimes \Z^2$ (see \cite[Theorem 1.9]{OT} and \cite{Mi}). Indeed, the Lie group $\R \ltimes \R^2$   is the group of matrices
 $$
 \left( \begin{array} {cccc} \cos (2 \pi t )& \sin(2 \pi t)& 0 & x\\
 -\sin (2 \pi t )& \cos(2 \pi t)& 0 & y\\
 0&0&1&t\\
 0&0&0&1 \end{array}
 \right)
 $$ and the lattice $\Gamma$  generated by $1$ in $\R$  and the standard lattice $\Z^2$, as in Example \ref{ablattice}.

By applying the Main Theorem    we have  that ${\mathcal A} (Ad_G (G)) = S^1 \ltimes \R^2$ and ${\mathcal A} (Ad_G (\Gamma)) = \R^2$.
Therefore  in this case $\tilde G \cong  \R^3 \subset S^1 \ltimes G$.  By applying the Main Theorem  we get that $\tilde G \cong \R^3$. Indeed, it is well known that $G / \Gamma$ is diffeomorphic to a torus.}
\end{ex}

The previous example $\R \ltimes \R^2$  is a special case of an {\em almost abelian}  Lie group. We recall that a Lie algebra ${\mathfrak g}$ is called almost abelian if it has an abelian  ideal of codimension $1$, i.e. if it can be represented as a semidirect sum $\R \ltimes {\mathfrak b}$, where ${\mathfrak b}$ is an abelian ideal of ${\mathfrak g}$.

In \cite{Gor} Gorbatsevich found a sufficient and necessary condition for which  given an almost abelian Lie group $G$ and a lattice $\Gamma$ of $G$ one has the equality ${\mathcal A} (Ad_G G)  = {\mathcal A} (Ad_G \Gamma)$. More precisely, he showed  
 the following

\begin{theorem} \cite[Theorem 4]{Gor} \label{ipirep} Let $G = \R \ltimes_{\varphi}  \R^n$ be a simply connected almost abelian Lie group and $\Gamma = \Z \ltimes \Z^n$
be a lattice in $G$.  Let $z$ be a generator for the subgroup $\Z \subset \Gamma$. The action of $z$ on $\R^n$ defines some matrix, i.e.,
 $z \in GL({\mathfrak g})$ and suppose that the one-parameter subgroup corresponding to the subgroup $A = \R$  is
${\mbox {exp}} (t á Z)$, where $Z$  is some matrix in the Lie algebra of derivations of the Lie algebra ${\mathfrak n}$ of the nilradical $N$ of $G$. Therefore  there is a unipotent
matrix $J$ such that $z = J á {\mbox {exp}}(Z)$.  

  Then ${\mathcal A}  (Ad_G G)={\mathcal A}  (Ad_G \Gamma)$   if and only if
the number $i \pi$  is not representable as a linear combination of the numbers  $\lambda_k$, $k = 1, \ldots, n$,   with rational coefficients, where $\{ 	\lambda_1, \ldots, \lambda_n \}$ is the spectrum of $Z$.
\end{theorem}

The idea for the proof of previous theorem is  that $$
\begin{array} {l}
{\mathcal A}( Ad_G (G)) = {\mathcal A} (Ad_G(A)) \ltimes Ad_G( N),\\[3pt]
  {\mathcal A}( Ad_G (\Gamma)) = {\mathcal A} (Ad_G(\Z)) \ltimes Ad_G( N),
\end{array}
$$
 since $Ad_G( N)$ is algebraic closed. Therefore,  finding  ${\mathcal A}( Ad_G (\Gamma))$  reduces to finding the
algebraic closure of the image of a cyclic subgroup (i.e. $\Z$  in the above decomposition for $\Gamma$).

In comparing ${\mathcal A} (Ad_G(A))$ and  ${\mathcal A} (Ad_G(\Z))$ the unipotent matrix $J$ plays no role, so one can suppose that $z = {\mbox {exp}}(Z)$. Therefore, we may assume  that
$A = {\mbox {exp}}(t Z)$ is the one-parameter subgroup generated by $Z$ and $\Z$ is the cyclic subgroup generated by $z = {\mbox {exp}}(Z)$.

One can use the decomposition   $z = z_s z_u$, where  $z_s$ (the semisimple part) is diagonal, i.e.  $z_s = {\mbox {diag}} ({\mbox {exp}} (\lambda_1, \ldots, \lambda_n))$ and $z_u$ is unipotent. The  semisimple part of the algebraic closure of $\Z$ consists of diagonal matrices and then it is contained in some algebraic torus. It is also clear that the algebraic closure of the  semisimple part $z_s$ is the intersection of the kernels of the diagonal matrices $\Pi_k e^{\lambda_k n_k}$.

In this way, Gorbatsevich showed that ${\mathcal A} (Ad_G(A)) ={\mathcal A} (Ad_G(\Z))$ if and only if $i \pi$ is representable as a linear combination of the numbers  $\lambda_k$, $k = 1, \ldots, n$,   with rational coefficients.

If $i \pi$ is not representable as a linear combination of the numbers  $\lambda_k$, $k = 1, \ldots, n$,   with rational coefficients, then one may apply Theorem \ref{Mostcoh} to compute the de Rham cohomology of the compact solvmanifold $G/\Gamma$ and one has that $H^*_{dR} (G/\Gamma) \cong H^*({\mathfrak g})$.

Otherwise, i.e. if $i \pi$  is   a rational  linear combination of the numbers  $\lambda_k$, $k = 1, \ldots, n$,   the only known result about the de Rham cohomology of $G/\Gamma$  is that  (see \cite[Proposition 4.7]{Bock}) $$b_1 (G/\Gamma) = n + 1 - {\mbox {rank}} (\varphi(1) - id).$$

By applying the Main Theorem  one obtains  a method to compute the de Rham cohomology of $G/\Gamma$ even in the latter  case.

It is not restrictive to suppose that ${\mathcal A}  (Ad_G \Gamma)$ is connected (otherwise one considers a finite index subgroup $\tilde \Gamma$ of $\Gamma$). Then ${\mathcal A}  (Ad_G G)= S^1  {\mathcal A}  (Ad_G \Gamma)$, or equivalently  ${\mathcal A} (Ad_G (A)) = S^1$,  and we may apply  the Main Theorem.

\begin{ex} \label{hyperell} {\bf (Hyperelliptic surface)} {\rm  Consider  the  solvable Lie group $G = \R \ltimes_{\varphi} (\C \times \R),$ where the action
$\varphi : \R  \to  Aut(\C \times \R)$  is defined by
$$
\varphi(t)(z, s) = (e^{i \eta   t} z, s),
$$
where $\eta = \pi, \frac {2}{3} \pi, \frac 12 \pi$ or $\frac 13 \pi$.

By \cite{Hasegawa} $G$ has seven isomorphism classes of  lattices $\Gamma = \Z \ltimes_{\varphi}  \Z^3$ where the action $\varphi: \Z \to Aut(\Z^3)$ is defined by a matrix $\varphi (1)$ which has   eigenvalues $1$, $e^{i \eta}$ and $e^{-i \eta}$.

Since $\varphi(1)$ has a pair of complex conjugate imaginary roots,  by applying Theorem \ref{ipirep}, in this case we have  ${\mathcal A} (Ad_G (G)) \neq {\mathcal A} (Ad_G (\Gamma))$. In this case ${\mathcal A} (Ad_G (\Gamma))$ is not connected, but  we have that $\Gamma$ 
 contains as a finite index subgroup $\tilde \Gamma \cong \Z^4$. In this way we get that the compact solvmanifold $G/\Gamma$ is a finite covering  of a torus.
 
Note that in this case $H^1_{dR} (G/\Gamma) \cong H^1 ({\mathfrak g})$ even if $G$ and $\Gamma$ do not satisfy the  Mostow condition. Indeed, $G$ has structure equations:
$$
\left \{ \begin{array} {l}
d e^1 = e^{24},\\
d e^2 = - e^{14},\\
d e^3 =0,\\
d e^4 =0
\end{array}
\right.
$$
and $H^1 ({\mathfrak g}) = {\mbox{span}} <e^3, e^4 >$.}

\end{ex}

\end{document}